\documentclass{amsart}
\usepackage{stmaryrd}
\usepackage{silence}
\usepackage{upgreek}
\usepackage{amssymb}
\usepackage{faktor}
\usepackage[all,cmtip]{xy}
\usepackage{amsmath} 
\usepackage{mathrsfs}
\usepackage{tikz}
\usepackage{tikz-cd}
\usepackage{enumitem}
\usepackage{mathtools}
\usepackage[mathcal]{euscript}
\usepackage{graphicx}
\usepackage{url}
\usepackage{hyperref}
\usepackage[T1]{fontenc}
\usepackage{subfigure}
\usepackage{epigraph}

\usetikzlibrary{shapes.geometric}
\usetikzlibrary{decorations.markings}
\usetikzlibrary{patterns,decorations.pathreplacing}

\usepackage{ragged2e}

\newcommand{\Map}{\mathrm{Map}}

\definecolor{coloryellow}{RGB}{240,228,66}
\definecolor{colorskyblue}{RGB}{86,180,233}
\definecolor{colorvermillion}{RGB}{213,94,0}

\theoremstyle{definition}
\newtheorem*{definition}{Definition}
\newtheorem*{problem}{Problem}

\theoremstyle{plain}
\newtheorem*{prototheorem}{Proto-theorem}
\newtheorem*{arnold relation}{Arnold relation}
\newtheorem*{theorem}{Theorem}

\usepackage{marginnote}
    \DeclareFontFamily{U}{wncy}{}
    \DeclareFontShape{U}{wncy}{m}{n}{<->wncyr10}{}
    \DeclareSymbolFont{mcy}{U}{wncy}{m}{n}
    \DeclareMathSymbol{\Sha}{\mathord}{mcy}{"58}

\newsavebox{\foobox}

\title{Lie algebras and the (co)homology of configuration spaces}
\author{Ben Knudsen}
\email{b.knudsen@northeastern.edu}
\address{Department of Mathematics, Northeastern University, Boston, MA 02115, USA}

\begin{document}

\begin{abstract}
We survey decades of research identifying the (co)homology of configuration spaces with Lie algebra (co)homology. The different routes to this one proto-theorem offer genuinely different explanations of its truth, and we attempt to convey some sense of the conceptual core of each perspective. We close with a list of problems.
\end{abstract}

\maketitle

\epigraph{\emph{I know noble accents\\
And lucid, inescapable rhythms;\\
But I know, too,\\
That the blackbird is involved\\
In what I know.\\
}
\hspace*{\fill}---Wallace Stevens}

\section{Introduction}

Having asked \emph{why}, a puzzle is receiving no answer; a mystery is receiving many. The mystery at hand concerns the (ordered) configuration space
\[F_k(X)=\{(x_1,\ldots, x_k)\in X^k: x_i\neq x_j \text{ if } i\neq j\}\] of $k$ particles in the background space $X$, together with its close cousin the unordered configuration space. The asking is prompted by results of the following type, which, as we will see, are legion.

\begin{prototheorem}
(Co)homology of configuration spaces is Lie algebra (co)homology.
\end{prototheorem}

This statement provokes many questions. What background space? What coefficients? What Lie algebra? What structure preserved by the isomorphism? But the fundamental question is \emph{why}. Why are there Lie algebras here at all?

The purpose of this paper is to present three genuinely different answers to this question. As these answers represent three major streams of the modern study of configuration spaces, the paper will be something of a whirlwind tour of the beautiful ideas and profound insights held in those streams. The only originality here, if any, is curatorial.

The reader may find it useful to know that the cohomology of a Lie algebra $\mathfrak{g}$ may be calculated from its Chevalley--Eilenberg complex, which is the (graded) symmetric algebra generated by a shift of $\mathfrak{g}$, equipped with the differential determined as a derivation by the bracket operation of $\mathfrak{g}$---see \cite[Ch. 7]{Weibel:IHA}, for example. The reader may also find it helpful eventually to have at least a passing familiarity with operads---see \cite{Vallette:AHO} and \cite{LodayVallette:AO} for gentle and systematic options, respectively. Later in the paper, we will make reference to the (contravariant) theory of operadic Koszul duality, first introduced in \cite{GinzburgKapranov:KDO}. Synthesizing and harmonizing the widely varying literature on this topic is beyond our scope, but see \cite{Fresse:KDOHPP,Heuts:KDCFG} for two more contemporary accounts.

\section{The Arnold relation}

In the study of configuration spaces of manifolds, as is so often the case, the first interesting example contains the generating principle. Recording the center of mass, distance from the center, and relative direction of particles determines a homeomorphism $F_2(\mathbb{R}^n)\cong\mathbb{R}^n \times \mathbb{R}_+\times S^{n-1}$; in particular, the \emph{Gauss map}
\begin{align*}
F_2(\mathbb{R}^n)&\xrightarrow{\gamma} S^{n-1}\\
(x_1,x_2)&\mapsto \frac{x_1-x_2}{|x_1-x_2|}
\end{align*}
is a homotopy equivalence. In general, we have a whole family of Gauss maps $\gamma_{ij}=\gamma\circ\pi_{ij}:F_k(\mathbb{R}^n)\to S^{n-1}$, where $\pi_{ij}$ is the projection onto the $i$th and $j$th coordinates, whence a family of cohomology classes $\alpha_{ij}=\gamma_{ij}^*(1)\in H^{n-1}(F_k(\mathbb{R}^n))$.

Using the foundational fact that the restricted coordinate projection $F_k(\mathbb{R}^n)\to F_{k-1}(\mathbb{R}^n)$ is a fiber bundle \cite[Thm. 1]{FadellNeuwirth:CS}, together with the Leray--Hirsch theorem, it is easy to show that these classes generate the cohomology ring. In fact, as shown in \cite{Arnold:CRCBG} for $n=2$ and \cite{Cohen:HCS}  in general, they freely generate up to a single non-obvious relation,\footnote{There are also the obvious relations $\alpha_{ji}=(-1)^{n}\alpha_{ij}$ and $\alpha_{ij}^2=0$ coming from $S^{n-1}$.} which bears a striking resemblance to the Jacobi identity of a Lie algebra---as we will see, for good reason.

\begin{arnold relation}
$\alpha_{ij}\alpha_{j\ell}+\alpha_{j\ell}\alpha_{\ell i}+\alpha_{\ell i}\alpha_{ij}=0.$
\end{arnold relation}

The remainder of this section will sketch three proofs of this relation, and, from this multiplicity, the main themes of the remainder will arise. We begin with the oldest proof and the origin of the eponym \cite{Arnold:CRCBG}.

\begin{proof}[Arnold's proof] Specializing to the case $n=2$,\footnote{Arnold's approach through de Rham cohomology admits a highly nontrivial extension to general $n$ involving the Fulton--MacPherson compactifications \cite{FultonMacPherson:CCS,Sinha:MTCCS} and graph complexes \cite{Kontsevich:OMDQ,LambrechtsVolic:FLNDO}.} we make the identification $\mathbb{R}^2=\mathbb{C}$. In de Rham cohomology, the class $\alpha_{ij}$ is represented by the $1$-form $\frac{1}{2\pi i}d\log(z_i-z_j)$. A little arithmetic shows that these $1$-forms themselves already satisfy the relation.
\end{proof}

Since the $1$-form appearing above measures the winding number around the diagonal where $z_i=z_j$, it is reasonable to think that this argument is premised on the combinatorics of diagonals. In contrast, our second argument will be premised on the combinatorics of projections \cite[\S 6]{Cohen:HCS}.

\begin{proof}[Cohen's proof] The group $H^{2(n-1)}(F_3(\mathbb{R}^n))$ is free Abelian of rank $2$ by the Leray--Hirsch argument alluded to above, so a relation of the form 
\[x\alpha_{12}\alpha_{23}+y\alpha_{23}\alpha_{31}+z\alpha_{31}\alpha_{12}=0\] must obtain. Application of the transpositions $\tau_{23}$ and $\tau_{12}$ to this relation forces the equalities $x=y$ and $x=z$. Since there is no torsion, and since the $\alpha_{ij}$ generate, the claim follows in this case, and the general case follows via pullback along the appropriate projection $F_k(\mathbb{R}^n)\to F_3(\mathbb{R}^n)$.
\end{proof}

The third argument will make essential use of manifold topology, in the form of Poincar\'{e} duality \cite{Sinha:HLDO}.

\begin{proof}[Sinha's proof] An equivalent formulation of the relation in question is the vanishing of the product $(\alpha_{ij}-\alpha_{i\ell})(\alpha_{j\ell}-\alpha_{ji}).$ Since $F_3(\mathbb{R}^n)$ is a manifold, this product is Poincar\'{e} dual to the transverse intersection of properly embedded submanifolds Poincar\'{e} dual to the factors (should such submanifolds exist). A little intersection theory shows that $\alpha_{ij}-\alpha_{i\ell}$ is Poincar\'{e} dual to the submanifold defined by requiring that $x_i$, $x_j$, and $x_\ell$ be collinear, with $x_i$ in between $x_j$ and $x_\ell$, which is disjoint from the submanifold with $x_j$ in between.
\end{proof}

As we will see, each of these perspectives on the Arnold relation offers its own explanation for the proto-theorem of the introduction, explanations premised on three genuinely different answers to a seemingly straightforward question: \emph{what is a Lie algebra?}

\section{Diagonals and partitions}

The combinatorics of diagonals is the combinatorics of partitions. More specifically, for every partition of the set $\{1,\ldots, k\}$ into blocks, there corresponds a generalized diagonal subspace of $M^k$, defined by requiring that the coordinates labeled by the elements of each block coincide. Accordingly, we may rephrase the calculation of Arnold and Cohen described in the previous section (additively) as the decomposition
\[H^*(F_k(\mathbb{R}^n))\cong \bigoplus_{\lambda \vdash k}\mathrm{Ind}_{\Sigma_\lambda}^{\Sigma_k}\bigotimes_{i=1}^{\ell(\lambda)}L_n(\lambda_i).\] Here $\lambda$ denotes an unordered partition of the number $k$ of length $\ell(\lambda)$, we make the abbreviation $L_n(m)=H^{(m-1)(n-1)}(F_m(\mathbb{R}^n))$, which is spanned as a $\Sigma_m$-module by the product class $\alpha_{12}\alpha_{23}\cdots\alpha_{m-1,m}$, and $\Sigma_\lambda$ is the product of wreath products of symmetric groups determined by $\lambda$---for example, we have $\Sigma_{(3,2,2)}=\Sigma_3\times \Sigma_2\wr\Sigma_2$, thought of as a subgroup of $\Sigma_7$.

There is a powerful and influential generalization of this formula to manifolds more interesting than Euclidean space, versions of which were articulated independently by \cite{CohenTaylor:CGFCCFSCCS,BenderskyGitler:CCFS, Kriz:ORHTCS, Totaro:CSAV}. The approach of Totaro is to analyze the stalks of the pushforward of the constant sheaf at the generalized diagonal subspaces described above, thereby showing that the Leray spectral sequence for the inclusion $F_k(M)\subseteq M^k$ has the following form.
\begin{theorem}[Totaro]\label{thm:totaro}
Let $M$ be an orientable manifold. There is a spectral sequence
\[E^{*,*}_2\cong \bigoplus_{\lambda \vdash k}\mathrm{Ind}_{\Sigma_\lambda}^{\Sigma_k}\bigotimes_{i=1}^{\ell(\lambda)}H^*(M^{\lambda_i};L_n(\lambda_i))\implies H^*(F_k(M)),\] in which the first nonzero differential is determined (as a derivation) by multiplication by the diagonal class of $M$ under the isomorphisms of $H^*(M^2; L_n(2))$ and  $H^*(M;L_n(1))^{\otimes 2}$ with $H^*(M^2)$.
\end{theorem}

As first observed by Getzler \cite{Getzler:RMHMCS,Getzler:HGSTSNLAASVS}, this $E_2$-page is the Chevalley--Eilenberg complex of a Lie algebra. In order to explain this observation, we recall that a \emph{symmetric sequence}\footnote{Otherwise known as a \emph{species} or $\mathrm{FB}$-\emph{module}.} is a sequence $\mathcal{X}=\{\mathcal{X}_k\}_{k\geq0}$ of objects (in some background category), where the weight\footnote{Otherwise known as \emph{arity.}} $k$ object $\mathcal{X}_k$ comes equipped with a $\Sigma_k$-action. Given a background monoidal structure, symmetric sequences may be tensored together via the recipe
\[(\mathcal{X}\otimes\mathcal{Y})_k=\bigoplus_{i+j=k}\mathrm{Ind}_{\Sigma_i\times\Sigma_j}^{\Sigma_k} \mathcal{X}_i\otimes\mathcal{X}_j.\] Following the current terminological trend, we will refer to algebraic structures in this setting via the (admittedly overloaded) epithet ``twisted''---thus, twisted algebra, twisted commutative algebra, and so on.\footnote{Such structures are discussed elsewhere in the literature (e.g., \cite{Fresse:KDOHPP}) using the language of left modules over operads.}

For us, the motivating examples of symmetric sequences are all functorial constructions derived from the symmetric sequence $F(X)_k=F_k(X)$ of spaces. In this language, we may reformulate the Arnold--Cohen calculation as the isomorphism 
\[H^*(F(\mathbb{R}^n))\cong\mathrm{Sym}(L_n)\] of symmetric sequences of graded modules.\footnote{In fact, as twisted commutative algebras---see Section \ref{section:projections}.} Here one can already see the germ of Getzler's result, which is a cochain level identification (over $\mathbb{R}$) with the Chevalley--Eilenberg complex of a the differential graded (twisted) Lie algebra obtained by tensoring the de Rham forms of $M$ with a free twisted Lie algebra.

More than a beautiful reformulation, this is a point of view with real teeth. First, it shows that the higher differentials in the Leray spectral sequence are essentially Massey products in $M$, systematizing the results of \cite{FelixThomas:CSMP}; in particular, the spectral sequence collapses for formal manifolds but not in general. Second, it is not difficult to show that the invariant part of this twisted Lie algebra is always formal, regardless of whether $M$ is, implying that the induced spectral sequence converging to the rational cohomology of unordered configuration spaces collapses \cite{Knudsen:BNSCSVFH}. The resulting expression in terms of Lie algebra cohomology unifies and extends all prior rational results and permits extensive computation \cite{BoedigheimerCohen:RCCSS,BoedigheimerCohenTaylor:OHCS,FelixThomas:RBNCS, Church:HSCSM, Drummond-ColeKnudsen:BNCSS}.

Getzler's work was later extended by Petersen to fabulous levels of generality: first, in allowing arbitrary coefficients for cohomology; second, in considering diagonal complements other than classical configuration spaces; and, third, in removing the restriction that the background space be a manifold \cite{Petersen:CGCS}. The price of this last, radical extension is working instead with compactly supported cohomology, with the old results recovered via Poincar\'{e} duality. Even at this level of generality, the key player remains the set of partitions, with its partial order by refinement.

Having taken seriously the perspective of diagonals and partitions, two conclusions now seem difficult to avoid. First, save for the intercession of Poincar\'{e} duality in the eleventh hour, it would appear that manifold topology plays essentially no role in the story. Second, in the form of the poset of partitions, it would seem that we have located the source of the incursion of Lie algebras, and hence a rationale for the proto-theorem of the introduction; indeed, a Lie algebra may be \emph{defined} as an algebra over an operad built from partition posets \cite{Fresse:KDOHPP}.

\section{Projections and commutativity}\label{section:projections}

The cohomology of a Lie algebra carries the structure of a commutative algebra; indeed, the differential in the Chevalley--Eilenberg complex is a derivation by definition. Thus, from the discussion above, we conclude that the symmetric sequence $H^*(F(M))$ carries the structure of a twisted commutative algebra, or TCA for short. Confronted with this (co)homological structure, it is natural for us to ask: can we locate a topological source for it?

In concrete terms, an algebra structure on the symmetric sequence $\mathcal{A}$ consists of the data of $\Sigma_i\times\Sigma_j$-equivariant maps $\mu_{ij}:\mathcal{A}_i\otimes \mathcal{A}_j\to \mathcal{A}_k$ for $i,j,k\geq0$ with $i+j=k$, subject to natural associativity and unitality axioms. Commutativity in this context is the requirement that the block permutation of $\{1,\ldots, i\}$ and $\{i+1,\ldots, k\}$ should intertwine $\mu_{ij}$ and $\mu_{ji}$. Having made our confusion specific,\footnote{In the words of Paul Goerss.} the desired topological source is now easily identified; indeed, for \emph{any} space $X$, the coordinate projections $F_k(X)\to F_i(X)\times F_j(X)$ equip the symmetric sequence $H^*(F(X))$ with the structure of a TCA.\footnote{At the level of unordered configuration spaces, this structure has a Poincar\'{e} dual description in terms of ``superposition'' of configurations, which provides yet another approach to the proto-theorem \cite{RandalWilliams:CSCM}.}

What can we learn from this extra structure? For one thing, in view of the equality $F_1(X)=X$, the unit element in $H^0(X)$ and the universal property of the free TCA produce a canonical map of TCAs of the form $\mathcal{S}\to H^*(F(X))$, where $\mathcal{S}$ denotes the free TCA on a single generator of weight $1$. In this way, we may regard the target as an $\mathcal{S}$-module.\footnote{The structure of an $\mathcal{S}$-module is equivalent to that of a functor from the category of finite sets and injections \cite{SamSnowden:ITCA}, denoted $\mathrm{FI}$ in the literature on stability phenomena.} The crucial fact about the TCA $\mathcal{S}$ is that, like its cousin the polynomial ring, it satisfies an analogue of the Hilbert basis theorem \cite{Snowden:SSEDM,ChurchEllenbergFarbNagpal:FIMONR,ChurchEllenbergFarb:FIMSRSG}.

\begin{theorem}[Church--Ellenberg--Farb--Nagpal--Snowden]
The twisted commutative algebra $\mathcal{S}$ is Noetherian.\footnote{The reader should not be fooled by the plausibility of this statement. Noetherianity results for TCAs are extremely difficult and, consequently, quite rare \cite{NagpalSamSnowden:NSDTTCA}.}
\end{theorem}

The relevance of this result for our purposes is that, as the reader is encouraged to verify as an exercise, the $E^2$-page in Totaro's theorem is finitely generated over $\mathcal{S}$ in each total degree, so Noetherianity implies that $H^*(F(M))$ is itself finitely generated in each degree. As explained in \cite{ChurchEllenbergFarb:FIMSRSG}, this fact places strong constraints on the behavior of the cohomology of $F_k(M)$ for large $k$. Known as \emph{representation stability}, this phenomenon implies that the multiplicities of irreducible symmetric group representations (working rationally) must eventually be constant; in particular, the Betti numbers of the unordered configuration spaces stabilize.\footnote{This account is ahistorical, as it was the study of homological stability, begun in \cite{McDuff:CSPNP} and continued in \cite{Church:HSCSM, RandalWilliams:HSUCS}, that led to the discovery of representation stability. The investigation of these and related stability phenomena has since become a thriving mathematical field of its own \cite{RandalWilliamsWahl:HSAG,GalatiusRandal-Williams:HSMSHDMII,Gadish:RSFLS,MillerWilson:HORSOCSM,GalatiusKupersRandalWilliams:ECMCG,AnDrummond-ColeKnudsen:ESHGBG,ProudfootRamos:CCG,KnudsenMillerTosteson:ESCS,BibbyGadish:GFANRSPOCS}.}

Returning to our main theme, we observe that our interpretation of $H^*(F(X))$ as a TCA renders unsurprising the appearance of Lie algebras and their cohomology in the story; indeed, according to \cite{Knudsen:PSTLA}, essentially \emph{every} TCA\footnote{Reduced, of finite type.} is quasi-isomorphic to a Chevalley--Eilenberg complex! This fact holds in such generality that it applies to the cohomology of a wide range of interconnected families of spaces, including Petersen's generalized configuration spaces as well as many others having no connection to diagonals or partitions.

We emerge from this contemplation of projections and commutativity with our understanding and our confusion augmented in equal measure. While our previous conviction that manifold topology is ancillary is bolstered, our belief in the primacy of partitions is called into question. With Lie algebras arising where no partition is in view, surely we must conclude that their true origin in our story lies in commutativity; indeed, a Lie algebra may be \emph{defined} as an algebra over the Koszul dual of the commutative operad \cite{GinzburgKapranov:KDO}.

\section{Poincar\'{e} duality and embeddings}

We come now to our third version of the narrative, which, far from downplaying the role of manifold topology, views it as central. The story begins with the observation that the space $F_k(\mathbb{R}^n)$ is homotopy equivalent to the subspace $\mathcal{E}_n(k)\subseteq\mathrm{Emb}(\sqcup_k \mathbb{R}^n,\mathbb{R}^n)$ of so-called rectilinear embeddings or \emph{little cubes}---viewing $\mathbb{R}^n$ as the interior of an $n$-dimensional cube, we allow only translation and scaling of sides \cite{May:GILS}. Since embeddings compose, the symmetric sequence $\mathcal{E}_n$ carries the structure of an operad, which is inherited by its homology $\mathrm{e}_n:=H_*(\mathcal{E}_n)\cong H_*(F(\mathbb{R}^n))$. 

\begin{theorem}[Cohen]
The operad $\mathrm{e}_n$ is isomorphic to the $n$-shifted Poisson operad.\footnote{This interpretation of Cohen's calculation is given in \cite[Thm. 1.6]{GetzlerJones:OHAIIDLS}.}
\end{theorem}

In other words, an $\mathrm{e}_n$-algebra has a commutative product and a Lie bracket of degree $n-1$ that is a derivation of the product. Plainly, then, this operad contains within it a shifted copy of the Lie operad, which, under the identification with the homology of configuration spaces, is nothing other than the dual of the symmetric sequence $L_n$ considered above. From this perspective, the Arnold relation is not merely reminiscent of the Jacobi identity; the two are equivalent.\footnote{See \cite{Sinha:HLDO} for a beautiful geometric description of the perfect pairing identifying the two relations.}

The discovery of $\mathcal{E}_n$ was motivated by the study of $n$-fold loop spaces, on which it acts. From this perspective, the configuration spaces of $\mathbb{R}^n$ together form a free algebra (up to homotopy), and this universal property leads to an approximation of $\Omega^n\Sigma^nX$ by the configuration space $F_X(\mathbb{R}^n)$ of particles labeled by the pointed space $X$, in which a particle labeled by the basepoint disappears. This relationship is globalized in the so-called \emph{scanning} or \emph{electric field map} 
\[F_X(M)\to \Gamma_c(E_X)\] of Segal \cite{Segal:CCT} and  McDuff \cite{McDuff:CSPNP}, whose target is the space of compactly supported sections of the $\Sigma^nX$-bundle associated to the frame bundle of $M$, where $O(n)$ acts on the suspension coordinates (when $M$ is parallelizable, this space of sections is simply $\Map_c(M,\Sigma^nX)$). 

\begin{theorem}[McDuff]
If $X$ is connected, then the scanning map is a weak equivalence.\footnote{See \cite{Boedigheimer:SSMS} for an exposition at this level of generality.}
\end{theorem}

It was recognized early on that the source of this map is a kind of homology theory for manifolds, being equipped with covariant functoriality (for open embeddings), monoidality, and an analogue of the long exact sequence of a pair. The target, on the other hand, is clearly a form of ``nonabelian'' compactly supported cohomology, so it is difficult to avoid thinking of this result as a form of Poincar\'{e} or Atiyah duality.\footnote{More recently, these ideas have been given systematic treatment in the theories of factorization homology \cite{AyalaFrancis:FHTM} and nonabelian Poincar\'{e} duality \cite{Salvatore:CSSL, Lurie:HA}.} Indeed, its proof follows a local-to-global blueprint entirely parallel to that of Poincar\'{e} duality.

Via the filtration of $F_X(M)$ by cardinality, McDuff's theorem permits another avenue of attack on the ordinary (unordered) configuration spaces of general manifolds \cite{BoedigheimerCohen:RCCSS,BoedigheimerCohenTaylor:OHCS,FelixThomas:RBNCS}. This last reference contains yet another version of the proto-theorem on Lie algebra cohomology; here it is obtained by applying rational homotopy theory \cite{Quillen:RHT,FelixHalperinThomas:RHT} to the target of the scanning map using results of Haefliger on rational models for mapping spaces \cite{Haefliger:RHSSNB}. In this way, we see the Chevalley--Eilenberg complex arising through a blending of Poincar\'{e} duality and Koszul duality,\footnote{See \cite{AyalaFrancis:PKD} for a systematic development of this idea.} the latter being the heart of rational homotopy theory.

This blending of dualities was first articulated at the level of operads in the following influential result \cite[Thm. 3.1]{GetzlerJones:OHAIIDLS}.

\begin{theorem}[Getzler--Jones]
In characteristic zero, the Koszul dual of $\mathrm{e}_n$ is $s^{-n}\mathrm{e}_n$.
\end{theorem}

Here we write $s^{-1}$ for the operadic desuspension, defined so that an $s^{-1}\mathcal{O}$-algebra structure on $X$ is simply an $\mathcal{O}$-algebra structure on the suspension of $X$.

As a consequence of this theorem, the obvious map of operads $\mathrm{e}_n\to \mathrm{e}_{n+1}$, induced by the standard inclusion between Euclidean spaces, gives rise to a map $s^{-n}\mathrm{e}_{n+1}\to s^{-n+1}\mathrm{e}_n$. In practical terms, the effect of these shifts is to place the top degree homology, which we earlier saw identified with the Lie operad, in degree $0$, independent of $n$. We now see that this manifestation of Poincar\'{e} duality is a perfectly viable explanation for the presence of Lie algebras in our study; indeed, it follows from Cohen's calculation that a Lie algebra may be \emph{defined} as an algebra over the inverse limit operad
\[\varprojlim\left(\,\cdots \to s^{-n}\mathrm{e}_{n+1}\to s^{-n+1}\mathrm{e}_n\to \cdots\to \mathrm{e}_1\right).\]

\section{The view from the sphere}

The self-duality theorem of Getzler--Jones was subsequently lifted to the level of operads in integral chain complexes by Fresse \cite{Fresse:KDEO}, then to the level of operads in spectra by Ching--Salvatore \cite{ChingSalvatore:KDTEO}, permitting the following definition.

\begin{definition}
The \emph{spectral Lie operad}\footnote{This definition is ahistorical---see \cite{Ching:BCTOGDI} for the original approach. Note that we adopt the grading convention in which the Lie bracket has degree $0$.} is the operad in spectra given by the inverse limit
\[\mathcal{L}=\varprojlim\left(\,\cdots \to s^{-n}\Sigma_+^\infty\mathcal{E}_{n+1}\to s^{-n+1}\Sigma_+^\infty\mathcal{E}_n\to \cdots\to \Sigma_+^\infty\mathcal{E}_1\right).\]
\end{definition}

Definitionally, then, there is a map of operads from the (shifted) spectral Lie operad to $\Sigma_+^\infty\mathcal{E}_n$, which results in an adjunction at the level of their categories of algebras, the left adjoint of which should be thought of as a higher enveloping algebra functor. The author studied these algebras in \cite{Knudsen:HEA}, in particular establishing their satisfaction of an analogue of the Poincar\'{e}--Birkhoff--Witt theorem.\footnote{Our account here is highly ahistorical. At the time of the writing of \cite{Knudsen:HEA}, spectral self-duality was still conjectural, and the author was forced to resort to a rather convoluted workaround inspired by the theory of chiral algebras \cite{BeilinsonDrinfeld:CA}, as interpreted in \cite{FrancisGaitsgory:CKD}. A proof of the equivalence of these two approaches to defining higher enveloping algebras has recently been announced in \cite{Antolin-CamarenaBrantnerHeuts:PBWTHA}.} Combining this result with the local-to-global philosophy of factorization homology \cite{AyalaFrancis:FHTM} yields the following result.\footnote{See the proof of \cite[Thm. C]{Knudsen:HEA}.}

\begin{theorem}[K]
Let $M$ be a smooth\footnote{For simplicity.} $n$-manifold of finite type. There is a weak equivalence \[
\Sigma_+^\infty F(M)\simeq \Sigma_+\mathrm{Bar}\left(\mathrm{id}, \mathcal{L}, \Map_c^{O(n)}\left(\mathrm{Fr}_M,\mathcal{L}( \Sigma^\infty S^{n-1})\right)\right),\]
of symmetric sequences of spectra, where $\mathrm{Bar}$ denotes the geometric realization of the simplicial two-sided bar construction, $\mathrm{Fr}_M$ the frame bundle of $M$,  and $\mathcal{L}$ the free spectral Lie algebra monad, and we view $S^{n-1}$ as a symmetric sequence concentrated in weight $1$.\footnote{In viewing the mapping spectrum as a spectral Lie algebra, we use that the $\infty$-category of spectral Lie algebra is cotensored over spaces.}
\end{theorem}

The bar construction being a spectral form of Lie algebra homology, this formula should be viewed as an articulation of the proto-theorem over the sphere spectrum, simultaneously recovering and expanding all prior (additive) results on the (co)homology of configuration spaces. By inspecting the formula, we may in particular conclude that the stable homotopy types of configuration spaces (ordered or unordered) of manifolds of fixed dimension are proper homotopy invariants.\footnote{See \cite{AouinaKlein:HICS} for prior partial results.}

The computational import of the theorem derives from the spectral sequence resulting from smashing with a (co)homology theory $E$ and filtering the bar construction by its skeleta. Converging to the $E$-(co)homology of configuration spaces,\footnote{Using a more general version of the theorem applying to labeled configuration spaces, and invoking McDuff's theorem, we equally obtain spectral sequences converging to the $E$-(co)homology of various mapping and section spaces---in particular, of iterated loop spaces.} the initial page of this spectral sequence is a form of Lie algebra (co)homology taking into account the extra information of power operations for spectral Lie algebras over $E$. 

In this way, each choice of $E$ gives rise to a different, typically quite intricate, computational problem, the solution to which lies in incorporating the earlier views of the Lie operad through partitions and Koszul duality. This study of power operations has been undertaken for ordinary mod $p$ homology \cite{AroneMahowald:GTIFUPHS,Kjaer:OPHFAOSLO,AntolinCamarena:MTHFSLA} and Morava $E$-theory \cite{Brantner:LTTSLA}, and the corresponding spectral sequence calculations begun in \cite{Zhang:QHSLAAMPHLCS,ChenZhang:MPHUCSS} and \cite{BrantnerHahnKnudsen:LTTCSI}, respectively, but there is much more to be done.

\section{Outlook}

We close with a selection of open problems. Although the author hopes one day to solve them, he would not mind if the reader got there first.

The first three are computational problems. That they remain unsolved may be surprising to the uninitiated (as it once was to the author).

\begin{problem}
Calculate the rational homology of the ordered configuration spaces of the torus.
\end{problem}

Even the stable homology is unknown.\footnote{See \cite{Pagaria:AGBNOCSEC} for partial results in this direction.} In fact, no stable multiplicity is known for any nontrivial irreducible representation of the symmetric group \cite[Prob. 3.5]{Farb:RS}.\footnote{The multiplicity of the trivial representation, stable and unstable, was calculated in \cite{Drummond-ColeKnudsen:BNCSS}.}

\begin{problem}
Calculate the mod $p$ homology of the unordered configuration spaces of a closed, orientable surface of positive genus, where $p$ is an odd prime.
\end{problem}

In contrast with the case of an open surface, treated in \cite{BianchiStavrouHCSSMOP}, almost nothing is known here.\footnote{See \cite{ChenZhang:MPHUCSS,BrantnerHahnKnudsen:LTTCSI} for some very partial results.}

\begin{problem}
Calculate the Morava $E$-theory and $K$-theory of $\Omega^kS^n$.\footnote{See \cite{Ravenel:WWSDKALSS} for some conjectures and \cite{Yamaguchi:MKDLSS,Langsetmo:KTLLOSA,Langsetmo:FSK,Tamaki:FIFSMKTO} for partial results.}
\end{problem}

Through McDuff's theorem, the spectral sequence described in the previous section provides one approach to this problem. Although the $E^2$-page has  a purely algebraic description \cite{BrantnerHahnKnudsen:LTTCSI}, the full computation seems difficult.

The final two problems are of a more conceptual leaning.

\begin{problem}
Is there an analogue of the Milnor--Moore theorem \cite{MilnorMoore:SHA} for higher enveloping algebras?
\end{problem}

In order to formulate the final problem, we recall that, up to homotopy, the ordered configuration spaces of a manifold carry the structure of a right module over the appropriate operad of little cubes.\footnote{We elide issues of parallelizability.}

\begin{problem}
Is the stable homotopy type of the operadic module of configuration spaces a proper homotopy invariant?\footnote{See \cite{Malin:SETOSCS} for partial results.}
\end{problem}

An affirmative answer to this problem would imply proper homotopy invariance of stable factorization homology \cite{AyalaFrancis:FHTM} and stable embedding calculus \cite{Weiss:EPVIT}.

\bibliographystyle{amsalpha}
\bibliography{references.bib}

\providecommand{\bysame}{\leavevmode\hbox to3em{\hrulefill}\thinspace}
\providecommand{\MR}{\relax\ifhmode\unskip\space\fi MR }
\providecommand{\MRhref}[2]{%
  \href{http://www.ams.org/mathscinet-getitem?mr=#1}{#2}
}
\providecommand{\href}[2]{#2}
\begin{thebibliography}{GKRW19}

\bibitem[AC20]{AntolinCamarena:MTHFSLA}
O.~Antolin-Camarena, \emph{The mod 2 homology of free spectral {Lie algebras}},
  Trans. Amer. Math. Soc. (2020).

\bibitem[ACBH]{Antolin-CamarenaBrantnerHeuts:PBWTHA}
O.~Antolin-Camarena, L.~Brantner, and G.~Heuts,
  \emph{{Poincar\'e--Birkhoff--Witt theorems in higher algebra}},
  arXiv:2501:03116.

\bibitem[ADCK20]{AnDrummond-ColeKnudsen:ESHGBG}
B.~H. An, G.~C. Drummond-Cole, and B.~Knudsen, \emph{Edge stabilization in the
  homology of graph braid groups}, Geom. Topol. \textbf{24} (2020).

\bibitem[AF15]{AyalaFrancis:FHTM}
D.~Ayala and J.~Francis, \emph{Factorization homology of topological
  manifolds}, J. Topol. \textbf{8} (2015), 1045--1084.

\bibitem[AF19]{AyalaFrancis:PKD}
\bysame, \emph{Poincar\'e/{K}oszul duality}, Comm. Math. Phys. \textbf{365}
  (2019).

\bibitem[AK04]{AouinaKlein:HICS}
M.~Aouina and J.~Klein, \emph{On the homotopy invariance of configuration
  spaces}, Algebr. Geom. Topol. \textbf{4} (2004), 813--827.

\bibitem[AM99]{AroneMahowald:GTIFUPHS}
G.~Arone and M.~Mahowald, \emph{The {Goodwillie} tower of the identity functor
  and the unstable periodic homotopy of spheres}, Invent. Math. \textbf{135}
  (1999), 743--788.

\bibitem[Arn69]{Arnold:CRCBG}
V.~I. Arnol'd, \emph{The cohomology ring of the colored braid group}, Math.
  Notes \textbf{5} (1969), no.~2, 138--140.

\bibitem[BC88]{BoedigheimerCohen:RCCSS}
C.-F. B{\"{o}}digheimer and F.~R. Cohen, \emph{Rational cohomology of
  configuration spaces of surfaces}, Algebraic Topology and Transformation
  Groups, Lecture Notes in Math., vol. 1361, Springer, 1988, pp.~7--13.

\bibitem[BCT89]{BoedigheimerCohenTaylor:OHCS}
C.-F. B{\"{o}}digheimer, F.~R. Cohen, and L.~R. Taylor, \emph{On the homology
  of configuration spaces}, Topology \textbf{28} (1989), 111--123.

\bibitem[BD69]{BeilinsonDrinfeld:CA}
A.~Beilinson and V.~Drinfeld, \emph{{Chiral Algebras}}, American Mathematical
  Society, 1969.

\bibitem[BG91]{BenderskyGitler:CCFS}
M.~Bendersky and S.~Gitler, \emph{The cohomology of certain function spaces},
  Trans. Amer. Math. Soc. \textbf{326} (1991).

\bibitem[BG23]{BibbyGadish:GFANRSPOCS}
C.~Bibby and N.~Gadish, \emph{A generating function approach to new
  representation stabillity phenomena in orbit configuration spaces}, Trans.
  Amer. Math. Soc. \textbf{10} (2023), 241--287.

\bibitem[BHK24]{BrantnerHahnKnudsen:LTTCSI}
L.~Brantner, J.~Hahn, and B.~Knudsen, \emph{The {Lubin--Tate} theory of
  configuration spaces: I}, J. Topol. (2024), to appear.

\bibitem[B{\"{o}}d87]{Boedigheimer:SSMS}
C.-F. B{\"{o}}digheimer, \emph{Stable splittings of mapping spaces}, Algebraic
  Topology, Lecture Notes in Math., vol. 1286, Springer, 1987.

\bibitem[Bra17]{Brantner:LTTSLA}
L.~Brantner, \emph{The {Lubin-Tate theory of spectral Lie} algebras}, Ph.D.
  thesis, Harvard University, 2017.

\bibitem[BS24]{BianchiStavrouHCSSMOP}
A.~Bianchi and A.~Stavrou, \emph{Homology of configuration spaces of surfaces
  modulo an odd prime}, J. Reine Angew. Math. (2024).

\bibitem[CEF15]{ChurchEllenbergFarb:FIMSRSG}
T.~Church, J.~S. Ellenberg, and B.~Farb, \emph{{FI}-modules and stability for
  representations of symmetric groups}, Duke Math. J. \textbf{164} (2015),
  no.~9, 1833--1910.

\bibitem[CEFN14]{ChurchEllenbergFarbNagpal:FIMONR}
T.~Church, J.~S. Ellenberg, B.~Farb, and R.~Nagpal, \emph{{FI}-modules over
  {Noetherian} rings}, Geom. Topol. (2014).

\bibitem[Chi05]{Ching:BCTOGDI}
M.~Ching, \emph{Bar constructions for topological operads and the goodwillie
  derivatives of the identity}, Geom. Topol. \textbf{9} (2005), 833--934.

\bibitem[Chu12]{Church:HSCSM}
T.~Church, \emph{Homological stability for configuration spaces of manifolds},
  Invent. Math. \textbf{188} (2012), no.~465-504.

\bibitem[Coh76]{Cohen:HCS}
F.~R. Cohen, \emph{The homology of {$\mathcal{C}_{n+1}$}-spaces, $n\geq0$},
  {The Homology of Iterated Loop Spaces}, Springer, 1976.

\bibitem[CS22]{ChingSalvatore:KDTEO}
M.~Ching and P.~Salvatore, \emph{Koszul duality for topological
  {$E_n$}-operads}, Proc. Lond. Math. Soc. (2022).

\bibitem[CT78]{CohenTaylor:CGFCCFSCCS}
F.~R. Cohen and L.~R. Taylor, \emph{{Computations of Gelfand--Fuks cohomology,
  the cohomology of function spaces, and the cohomology of configuration
  spaces}}, {Geometry Applications of Homotopy Theory I}, Lecture Notes in
  Math., vol. 657, Springer, 1978.

\bibitem[CZ22]{ChenZhang:MPHUCSS}
M.~Chen and A.~Zhang, \emph{{Mod $p$ homology of unordered configuration spaces
  of surfaces}}, Proc. Amer. Math. Soc (2022), to appear.

\bibitem[DCK17]{Drummond-ColeKnudsen:BNCSS}
G.~C. Drummond-Cole and B.~Knudsen, \emph{Betti numbers of configuration spaces
  of surfaces}, J. London Math. Soc. \textbf{96} (2017), no.~2, 367--393.

\bibitem[Far]{Farb:RS}
B.~Farb, \emph{Representation stability}, Contribution to the proceedings of
  the ICM 2014, Seoul, arXiv:1404.4065.

\bibitem[FG12]{FrancisGaitsgory:CKD}
J.~Francis and D.~Gaitsgory, \emph{Chiral {K}oszul duality}, Selecta Math.
  \textbf{18} (2012), no.~1, 27--87.

\bibitem[FHT00]{FelixHalperinThomas:RHT}
Y.~F{\'{e}}lix, S.~Halperin, and J.-C. Thomas, \emph{Rational homotopy theory},
  Springer, 2000.

\bibitem[FM94]{FultonMacPherson:CCS}
W.~Fulton and R.~MacPherson, \emph{A compactification of configuration spaces},
  Ann. of Math. (2) \textbf{139} (1994), no.~1, 183--225.

\bibitem[FN62]{FadellNeuwirth:CS}
E.~Fadell and L.~Neuwirth, \emph{Configuration spaces}, Math. Scand.
  \textbf{10} (1962), 111--118.

\bibitem[Fre04]{Fresse:KDOHPP}
B.~Fresse, \emph{Koszul duality of operads and homology of partition posets},
  Homotopy theory: relations with algebraic geometry, group cohomology, and
  algebraic {K}-theory, Contemp. Math., vol. 346, Amer. Math. Soc., 2004.

\bibitem[Fre10]{Fresse:KDEO}
\bysame, \emph{Koszul duality of {$E_n$}-operads}, Selecta Math. \textbf{17}
  (2010), 363--434.

\bibitem[FT00]{FelixThomas:RBNCS}
Y.~F\'elix and J.-C. Thomas, \emph{Rational {B}etti numbers of configuration
  spaces}, Topology Appl. \textbf{102} (2000), 139--149.

\bibitem[FT04]{FelixThomas:CSMP}
\bysame, \emph{Configuration spaces and {Massey} products}, Int. Math. Res.
  Not. \textbf{33} (2004), 1685--1702.

\bibitem[Gad17]{Gadish:RSFLS}
N.~Gadish, \emph{Representation stability for families of linear subspace
  arrangements}, Adv. Math. (2017), 341--377.

\bibitem[Get]{Getzler:HGSTSNLAASVS}
E.~Getzler, \emph{The homology groups of some two-step nilpotent lie algebras
  associated to symplectic vector spaces}, arXiv:9903147.

\bibitem[Get99]{Getzler:RMHMCS}
\bysame, \emph{Resolving mixed {H}odge modules on configuration spaces}, Duke
  Math. J. \textbf{96} (1999), no.~1, 175--203.

\bibitem[GJ]{GetzlerJones:OHAIIDLS}
E.~Getzler and J.~D.~S. Jones, \emph{Operads, homotopy algebra and iterated
  integrals for double loop spaces}, arXiv:hep-th/9403055v1.

\bibitem[GK94]{GinzburgKapranov:KDO}
V.~Ginzburg and M.~M. Kapranov, \emph{Koszul duality for operads}, Duke Math.
  J. \textbf{76} (1994), 203--272.

\bibitem[GKRW19]{GalatiusKupersRandalWilliams:ECMCG}
S.~Galatius, A.~Kupers, and O.~Randal-Williams, \emph{{$E_2$}-cells and mapping
  class groups}, Publ. Math. Inst. Hautes \'{E}tudes Sci. \textbf{130} (2019),
  1--61.

\bibitem[GRW17]{GalatiusRandal-Williams:HSMSHDMII}
S.~Galatius and O.~Randal-Williams, \emph{Homological stability for moduli
  spaces of high dimensional manifolds. ii}, Ann. of Math. (2) \textbf{186}
  (2017).

\bibitem[Hae82]{Haefliger:RHSSNB}
A.~Haefliger, \emph{Rational homotopy of the space of sections of a nilpotent
  bundle}, Trans. Amer. Math. Soc. \textbf{273} (1982), 609--620.

\bibitem[Heu]{Heuts:KDCFG}
G.~Heuts, \emph{{Koszul duality and a conjecture of Francis--Gaitsgory}},
  arXiv:2408.06173.

\bibitem[Kja18]{Kjaer:OPHFAOSLO}
J.~Kjaer, \emph{On the odd primary homology of free algebras over the spectral
  {Lie} operad}, J. Homotopy Relat. Str. \textbf{13} (2018), 581--597.

\bibitem[KMT22]{KnudsenMillerTosteson:ESCS}
B.~Knudsen, J.~Miller, and P.~Tosteson, \emph{Extremal stability for
  configuration spaces}, Math. Ann. (2022).

\bibitem[Knu17]{Knudsen:BNSCSVFH}
B.~Knudsen, \emph{Betti numbers and stability for configuration spaces via
  factorization homology}, Alg. Geom. Topol. \textbf{17} (2017), no.~5,
  3137--3187.

\bibitem[Knu18]{Knudsen:HEA}
\bysame, \emph{Higher enveloping algebras}, Geom. Topol. \textbf{22} (2018),
  no.~7.

\bibitem[Knu22]{Knudsen:PSTLA}
\bysame, \emph{Projection spaces and twisted {Lie} algebras}, Contemp. Math.
  (2022), To appear.

\bibitem[Kon99]{Kontsevich:OMDQ}
M.~Kontsevich, \emph{Operads and motives in deformation quantization}, Lett.
  Math. Phys. \textbf{48} (1999), 35--72.

\bibitem[Kri94]{Kriz:ORHTCS}
I.~Kriz, \emph{On the rational homotopy type of configuration spaces}, Ann. of
  Math. (2) \textbf{139} (1994).

\bibitem[Lan93]{Langsetmo:KTLLOSA}
L.~Langsetmo, \emph{The {$K$}-theory localization of loops on an odd sphere and
  applications}, Topology \textbf{32} (1993).

\bibitem[Lan96]{Langsetmo:FSK}
\bysame, \emph{Further structure in {$K(1)_*\Omega^kS^{2n+1}$}}, {Algebraic
  Topology: New Trends in Localization and Periodicity}, Progress in
  Mathematics, vol. 136, Springer, 1996.

\bibitem[Lur]{Lurie:HA}
J.~Lurie, \emph{Higher algebra},
  \url{http://www.math.harvard.edu/~lurie/papers/higheralgebra.pdf}.

\bibitem[LV12]{LodayVallette:AO}
J.-L. Loday and B.~Vallette, \emph{Algebraic operads}, Grundlehren Math. Wiss.,
  vol. 346, Springer Verlag, 2012.

\bibitem[LV13]{LambrechtsVolic:FLNDO}
P.~Lambrechts and I.~Voli\'{c}, \emph{Formality of the little {$N$}-disks
  operad}, Memoirs of the American Mathematical Society, Amer. Math. Soc.,
  2013.

\bibitem[Mal23]{Malin:SETOSCS}
C.~Malin, \emph{The stable embedding tower and operadic structures on
  configuration spaces}, Homol. Homot. Appl. (2023).

\bibitem[May72]{May:GILS}
J.~P. May, \emph{{The Geometry of Iterated Loop Spaces}}, Lecture Notes in
  Math., vol. 271, Springer-Verlag, Berlin, Germany, 1972.

\bibitem[McD75]{McDuff:CSPNP}
D.~McDuff, \emph{Configuration spaces of positive and negative particles},
  Topology \textbf{14} (1975), 91--107.

\bibitem[MM65]{MilnorMoore:SHA}
J.~Milnor and J.~C. Moore, \emph{On the structure of {Hopf} algebras}, Ann. of
  Math. (2) \textbf{81} (1965), no.~2, 211--264.

\bibitem[MW19]{MillerWilson:HORSOCSM}
J.~Miller and J.~Wilson, \emph{Higher order representation stability and
  ordered configuration spaces of manifolds}, Geom. Topol. \textbf{23} (2019).

\bibitem[NSS16]{NagpalSamSnowden:NSDTTCA}
R.~Nagpal, S.~V. Sam, and A.~Snowden, \emph{Noetherianity of some degree two
  twisted commutative algebras}, Selecta Math. \textbf{22} (2016), no.~2,
  913--937.

\bibitem[Pag22]{Pagaria:AGBNOCSEC}
R.~Pagaria, \emph{Asymptotic growth of {Betti} numbers of ordered configuration
  spaces of an elliptic curve}, Eur. J. Math. \textbf{8} (2022), 427--445.

\bibitem[Pet20]{Petersen:CGCS}
D.~Petersen, \emph{Cohomology of generalized configuration spaces}, Compos.
  Math. \textbf{156} (2020), no.~2, 251--298.

\bibitem[PR22]{ProudfootRamos:CCG}
N.~Proudfoot and E.~Ramos, \emph{The contraction category of graphs},
  Represent. Theory (2022).

\bibitem[Qui69]{Quillen:RHT}
D.~Quillen, \emph{Rational homotopy theory}, Ann. of Math. (2) \textbf{90}
  (1969), 205--295.

\bibitem[Rav98]{Ravenel:WWSDKALSS}
D.~Ravenel, \emph{What we still don't know about loop spaces of spheres},
  Homotopy theory via algebraic geometry and group representations,
  Contemporary Mathematics, vol. 220, Amer. Math. Soc., 1998.

\bibitem[RW13]{RandalWilliams:HSUCS}
O.~Randal-Williams, \emph{Homological stability for unordered configuration
  spaces}, Q. J. Math. \textbf{64} (2013), 303--326.

\bibitem[RW24]{RandalWilliams:CSCM}
\bysame, \emph{Configuration spaces as commutative monoids}, Bull. London Math.
  Soc. (2024).

\bibitem[RWW17]{RandalWilliamsWahl:HSAG}
O.~Randal-Williams and N.~Wahl, \emph{Homological stability for automorphism
  groups}, Adv. Math. (2017), no.~534--626.

\bibitem[Sal01]{Salvatore:CSSL}
P.~Salvatore, \emph{Configuration spaces with summable labels}, Cohomological
  Methods in Homotopy Theory, Progr. Math., vol. 196, Birkh\"auser, 2001.

\bibitem[Seg74]{Segal:CCT}
G.~Segal, \emph{Categories and cohomology theories}, Topology \textbf{13}
  (1974), no.~3, 293--312.

\bibitem[Sin]{Sinha:HLDO}
D.~Sinha, \emph{The homology of the little disks operad}, arXiv:0610236.

\bibitem[Sin04]{Sinha:MTCCS}
\bysame, \emph{Manifold-theoretic compactifications of configuration spaces},
  Selecta Math. \textbf{10} (2004), no.~3, 391--428.

\bibitem[Sno13]{Snowden:SSEDM}
A.~Snowden, \emph{Syzygies of {Segre} embeddings and {$\Delta$}-modules}, Duke
  Math. J. \textbf{162} (2013), 225--277.

\bibitem[SS]{SamSnowden:ITCA}
S.~V. Sam and A.~Snowden, \emph{Introduction to twisted commutative algebras},
  arXvi:1209.5122.

\bibitem[Tam02]{Tamaki:FIFSMKTO}
D.~Tamaki, \emph{The fiber of iterated freudenthal suspension and {Morava
  $K$-theory} of {$\Omega^kS^{2\ell+1}$}}, Contemp. Math. \textbf{293} (2002).

\bibitem[Tot96]{Totaro:CSAV}
B.~Totaro, \emph{Configuration spaces of algebraic varieties}, Topology
  \textbf{35} (1996), no.~4, 1057--1067.

\bibitem[Val14]{Vallette:AHO}
B.~Vallette, \emph{Algebra {$+$} homotopy {$=$} operad}, Symplectic, Poisson
  and Noncommutative Geometry, vol.~62, MSRI Publications, 2014.

\bibitem[Wei94]{Weibel:IHA}
C.~Weibel, \emph{{An Introduction to Homological Algebra}}, Cambridge
  University Press, 1994.

\bibitem[Wei99]{Weiss:EPVIT}
M.~Weiss, \emph{Embeddings from the point of view of immersion theory: part
  {I}}, Geom. Topol. \textbf{3} (1999), 67--101.

\bibitem[Yam88]{Yamaguchi:MKDLSS}
A.~Yamaguchi, \emph{Morava {$K$}-theory of iterated loop spaces of spheres},
  Math. Z. \textbf{199} (1988).

\bibitem[Zha21]{Zhang:QHSLAAMPHLCS}
A.~Zhang, \emph{{Quillen homology of spectral Lie algebras with application to
  mod $p$ homology of labeled configuration spaces}}, Alg. Geom. Topol. (2021),
  to appear.

\end{thebibliography}

\end{document}